\documentclass{amsart}
\usepackage{amsmath,amsthm,amssymb,IMjournal}
\usepackage{amsfonts}
\usepackage{exscale}
\usepackage{graphicx}
\usepackage{multirow}
\usepackage{url}
\usepackage[bookmarks=false]{hyperref}
\begin{document}
\newtheorem{theorem}{Theorem}[section]
\newtheorem{assumption}[theorem]{Assumption}
\newtheorem{algorithm}[theorem]{Algorithm}
\newtheorem{lemma}[theorem]{Lemma}
\newtheorem{remark}[theorem]{Remark}

\numberwithin{equation}{section}

\title{On Surrogate Learning for Linear Stability Assessment of Navier\textendash
Stokes Equations with Stochastic Viscosity \\
\centering
\small\emph{Dedicated to Professor Pavel Burda on the occasion of his 75th birthday.}
}

\author{\|Bed\v{r}ich |Soused\'{\i}k|, Baltimore, MD, USA \\
             \|Howard |Elman|, College Park, MD, USA \\
             \|Kookjin |Lee|, Tempe, AZ, USA \\
             \|Randy |Price|, Fairfax, VA, USA}

\rec {March 3, 2022}

\abstract
We study linear stability of solutions to the Navier\textendash Stokes equations with
stochastic viscosity. Specifically, we assume that the viscosity is given in
the form of a~stochastic expansion. Stability analysis requires a solution of
the steady-state Navier\textendash Stokes equation and then leads to a generalized
eigenvalue problem, from which we wish to characterize the real part of the
rightmost eigenvalue. While this can be achieved by Monte Carlo simulation,
due to its computational cost we study three surrogates based on generalized
polynomial chaos, Gaussian process regression and a shallow neural network.
The results of linear stability analysis assessment obtained by the surrogates
are compared to that of Monte Carlo simulation using a set of numerical experiments.
\endabstract



\thanks
The work was supported by the U.S. Department of Energy Office of Advanced
Scientific Computing Research, Applied Mathematics program under award DE-SC0009301
and by the U.S. National Science Foundation under awards DMS1819115 and DMS1913201.
Bed\v{r}ich Soused\'{i}k would like to thank Professors Pavel Burda and Jaroslav Novotn\'{y}
for many inspiring discussions about bifurcations and Navier\textendash Stokes equations that inspired this work.
Part of the work was completed while Randy Price was a graduate student at the University of Maryland, Baltimore County. 
\endthanks

\section{Introduction}

Models of mathematical physics are typically based on partial differential
equations and they are often solved numerically using finite element methods.
The models use parameters as input data, although exact parameter values are
often not known and they are modeled using random variables. This approach
leads to so-called partial differential equations with uncertain data: given stochastic
parameters, we wish to characterize their stochastic solutions. Probably the
most popular method for solution of these problems is Monte Carlo simulation,
which is 
based on sampling: samples of input parameters give a set of
independent deterministic problems, which are solved, and then the statistical
moments of solution are obtained from ensemble averaging.
This method is known to be slow (with errors for $n$ samples behaving like~$n^{-1/2}$),
and since each sample
requires solution of the full model, its computational costs will be high.
Significant effort has been devoted to design
of computationally cheaper alternatives to the full model called
\emph{surrogates} in order to decrease the overall computational cost.
Arguably the most popular surrogate types are based on generalized polynomial
chaos~(gPC) in the engineering community~\cite{Ghanem-1991-SFE,Xiu-2002-WAP}, 
and Gaussian process~(GP)\ regression in the statistics 
community~\cite{Rasmussen-2005-GPM,STein-1999-ISD}. 

Our focus is on linear stability analysis of parameterized dynamical systems.
A steady solution $u$ is {\em stable} if with a small perturbation of $u$, used
as initial data in a transient simulation, the simulation reverts to $u$; otherwise it is {\em unstable}.
This is of fundamental importance in studying dynamics, since unstable solutions 
may lead to turbulent flows or other inexplicable dynamic 
behavior~\cite{Cliffe-Spence-Taverner,Schmid-Henningson}.
{\em Linear stability analysis} entails computing the rightmost eigenvalue of the Jacobian matrix at $u$;
if this eigenvalue has positive real part, then $u$ is unstable.
In this study, we explore this issue using the parameterized Navier\textendash Stokes equations,
This is a challenging task, because it entails solution of a nonlinear PDE close 
to a bifurcation point 
followed by a solution of a nonsymmetric eigenvalue problem. 
Since it is also computationally intensive, we wish to find a less expensive surrogate. 
The Navier\textendash Stokes equations with stochastic viscosity were studied, e.g., 
by~\cite{Lee-2019-LRS,Powell-2012-PSS,Sousedik-2016-SGM} and techniques
based on gPC for parameterized eigenvalue problems were studied, e.g., 
by~\cite{Andreev-2012-STA,Benner-2019-LRI,Lee-2018-IMS}. 
A stochastic collocation method for linear stability analysis was studied in~\cite{Elman-2018-CME}.
Most recently, an algorithm for solving nonsymmetric eigenvalue problems with uncertain data using an embedded (intrusive) stochastic Galerkin method, 
and the same application as in the present study, was proposed in~\cite{Sousedik-2022-SGM}.

Specifically, we design and compare several surrogates. There is only a handful
of studies comparing gPC approaches\ and GP\ regression see, 
e.g.,~\cite{O'Hagan-2013-PCT,Owen-2017-CSB,Yan-2018-GPP}, one of our goals is to
contribute to the discussion with this particularly challenging problem. 
For the construction of the gPC\ surrogate we
use the stochastic collocation method, and in particular the variant based on
the pseudospectral (nonintrusive) stochastic Galerkin method, 
see~\cite{Babuska-2010-SCM,Xiu-2010-NMS}, and for the GP\ surrogate we use the
\textsc{Matlab} function \texttt{fitrgp}. We note that it seems quite common
to use software packages for GP\ regression, and very different results
among the packages have been reported~\cite{Erickson-2018-CGP}. Therefore,
in our numerical experiments we compare both gPC\ and GP\ surrogates with
results obtained from Monte Carlo simulation. Finally, following recent trends
in using neural networks for solving PDE-based models
see, e.g.,~\cite{Peng-2020-MMM,Sirignano-2018-DGM}, 
we study a surrogate based on a shallow neural network. We compare the performance of the
surrogates using two benchmark problems, and we also compare the results with
those of Monte Carlo simulation.

The paper is organized as follows. In Section~\ref{sec:NS} we recall the
Navier\textendash Stokes equations and the finite element discretization, in
Section~\ref{sec:stability} we discuss the linear stability of the model, in
Section~\ref{sec:stochastic} we formulate the Navier\textendash Stokes equations with
stochastic viscosity and introduce the surrogates, in
Section~\ref{sec:numerical} we present results of numerical experiments, and
in Section~\ref{sec:conclusion} we summarize our work.

\section{Steady-state Navier\textendash Stokes equations}

\label{sec:NS}We begin by defining the model and notation for the
deterministic steady-state Navier\textendash Stokes equations,
following~\cite{Elman-2014-FEF}. We wish to find velocity$~\vec{u}$ and
pressure$~p$ such that
\begin{align}
-\nu\nabla^{2}\vec{u}+\left(  \vec{u}\cdot\nabla\right)  \vec{u}+\nabla p  &
=\vec{f},\label{eq:NS-1}\\
\nabla\cdot\vec{u}  &  =0, \label{eq:NS-2}%
\end{align}
in a spatial domain$~D$, satisfying boundary conditions
\begin{equation} \label{Dir-Neu-bc}
\vec{u}   =\vec{g}\ \text{ on }\Gamma_{\text{Dir}}, \qquad
\nu\nabla\vec{u}\cdot\vec{n}-p\vec{n}   =\vec{0} \text{ on }%
\Gamma_{\text{Neu}}, 
\end{equation}
where$~\partial D=\overline{\Gamma}_{\text{Dir}}\cup\overline{\Gamma
}_{\text{Neu}}$, $\vec{n}$ denotes the normal vector, 
$\nu$ denotes the kinematic viscosity and $\vec{f}$ is a vector of external forces, 
and we assume sufficient
regularity of the data.
Properties of the flow are usually characterized by the Reynolds number
\begin{equation}
Re=\frac{UL}{\nu}, \label{eq:Re}%
\end{equation}
where $U$ is a characteristic velocity and $L$ a characteristic length. 

In the mixed variational formulation of~(\ref{eq:NS-1})--(\ref{eq:NS-2}) we
wish to find $\left(  \vec{u},p\right)  \in\left(  V_{E},Q_{D}\right)  $ such
that
\begin{align}
\int_{D}  \nu\, \nabla\vec{u}:\nabla\vec{v}+\int_{D}\left(  \vec{u}\cdot\nabla
\vec{u}\right)  \cdot\vec{v}-\int_{D}p\left(  \nabla\cdot\vec{v}\right)   &
=\int_{D}\vec{f}\cdot\vec{v},\quad\forall\vec{v}\in V_{D}%
,\label{eq:NS-variational-1}\\
\int_{D}q\left(  \nabla\cdot\vec{u}\right)   &  =0,\quad\forall q\in Q_{D},
\label{eq:NS-variational-2}%
\end{align}
where $\left(  V_{D},Q_{D}\right)  $ is a pair of spaces satisfying an inf-sup
condition and~$V_{E}$ is an extension of~$V_{D}$ containing velocity vectors
that satisfy the Dirichlet boundary conditions~\cite{Girault-1986-FEM}.

Let $c(\vec{z};\vec{u},\vec{v})\equiv\int_{D}\left(  \vec{z}\cdot\nabla\vec
{u}\right)  \cdot\vec{v}$. Because the problem~(\ref{eq:NS-variational-1}%
)--(\ref{eq:NS-variational-2}) is nonlinear, it is solved using a
linearization scheme in the form of Newton or Picard iteration, derived as
follows.\footnote{This gives direct computation of the steady solution.
It is also possible to find such solutions by integrating to steady state; 
see, for example~\cite{Akervik-2006-SSN,Loiseau-2019-TKM}.}
Consider a solution $\left(  \vec{u},p\right)  $
of~(\ref{eq:NS-variational-1})--(\ref{eq:NS-variational-2}) to be given as
$\vec{u}=\vec{u}^{n}+\delta\vec{u}^{n}$ and $p=p^{n}+\delta p^{n}$.
Substituting into~(\ref{eq:NS-variational-1})--(\ref{eq:NS-variational-2}) and
neglecting the quadratic term $c(\delta\vec{u}^{n};\delta\vec{u}^{n},\vec{v})$
gives
\begin{align}
\int_{D} \nu\, \nabla\delta\vec{u}^{n}:\nabla\vec{v}+c(\delta\vec{u}^{n};\vec
{u}^{n},\vec{v})+c(\vec{u}^{n};\delta\vec{u}^{n},\vec{v})-\int_{D}\delta
p^{n}\left(  \nabla\cdot\vec{v}\right)   &  =R^{n}\left(  \vec{v}\right)
,\label{eq:Newton-1}\\
\int_{D}q\left(  \nabla\cdot\delta\vec{u}^{n}\right)   &  =r^{n}\left(
q\right)  , \label{eq:Newton-2}%
\end{align}
where
\begin{align}
R^{n}\left(  \vec{v}\right)   &  =\int_{D}\vec{f}\cdot\vec{v}- \int_{D} \nu\, \nabla\vec{u}^{n}:\nabla\vec{v}-c(\vec{u}^{n};\vec{u}^{n},\vec{v}%
)+\int_{D}p^{n}\left(  \nabla\cdot\vec{v}\right)  ,\label{eq:residual-1}\\
r^{n}\left(  q\right)   &  =-\int_{D}q\left(  \nabla\cdot\vec{u}^{n}\right)  .
\label{eq:residual-2}%
\end{align}
Step~$n$\ of the \emph{Newton iteration} obtains $\left(  \delta\vec{u}%
^{n},\delta p^{n}\right)  $ from~(\ref{eq:Newton-1})--(\ref{eq:Newton-2}) and
updates the solution as
\begin{equation} \label{eq:Newton-increment}
\vec{u}^{n+1}    =\vec{u}^{n}+\delta\vec{u}^{n}, \qquad
p^{n+1}   =p^{n}+\delta p^{n}.
\end{equation}
Step~$n$\ of the \emph{Picard iteration} omits the term $c(\delta\vec{u}%
^{n};\vec{u}^{n},\vec{v})$ in (\ref{eq:Newton-1}), giving
\begin{align}
\int_{D} \nu\, \nabla\delta\vec{u}^{n}:\nabla\vec{v}+c(\vec{u}^{n};\delta\vec
{u}^{n},\vec{v})-\int_{D}\delta p^{n}\left(  \nabla\cdot\vec{v}\right)   &
=R^{n}\left(  \vec{v}\right)  ,\label{eq:Picard-1}\\
\int_{D}q\left(  \nabla\cdot\delta\vec{u}^{n}\right)   &  =r^{n}\left(
q\right)  . \label{eq:Picard-2}%
\end{align}

Next, let us consider the discretization of~(\ref{eq:NS-1})--(\ref{eq:NS-2})
by a div-stable mixed finite element method,
and let the bases for the velocity and pressure spaces be denoted $\left\{
\phi_{i}\right\}  _{i=1}^{n_{u}}$\ and $\left\{  \varphi_{i}\right\}
_{i=1}^{n_{p}}$, respectively, $n_{u}>n_{p}$, and let us denote by
$n_{x}=n_{u}+n_{p}$ the number of velocity and pressure degrees of freedom. In
matrix terminology, each nonlinear iteration entails solving a linear system
\begin{equation}
\left[
\begin{array}
[c]{cc}%
\mathbf{F}^{n} & \mathbf{B}^{T}\\
\mathbf{B} & \mathbf{0}%
\end{array}
\right]  \left[
\begin{array}
[c]{c}%
\delta\mathbf{u}^{n}\\
\delta\mathbf{p}^{n}%
\end{array}
\right]  =\left[
\begin{array}
[c]{c}%
\mathbf{R}^{n}\\
\mathbf{r}^{n}%
\end{array}
\right]  , \label{eq:Newton}%
\end{equation}
which corresponds to (\ref{eq:Newton-1})--(\ref{eq:Newton-2}), followed by an
update of the solution
\begin{equation} 
\label{eq:updata}
\mathbf{u}^{n+1}   =\mathbf{u}^{n}+\delta\mathbf{u}^{n}, \qquad
\mathbf{p}^{n+1}    =\mathbf{p}^{n}+\delta\mathbf{p}^{n}. 
\end{equation} 
For Newton's method, $\mathbf{F}^{n}$ is the (nonsymmetric)\ Jacobian matrix,
a sum of the vector-Laplacian matrix~$\mathbf{A}$, the vector-convection
matrix~$\mathbf{N}^{n}$, and the Newton derivative matrix~$\mathbf{W}^{n}$,
\begin{equation}
\mathbf{F}^{n}=\mathbf{A}+\mathbf{N}^{n}+\mathbf{W}^{n},
\label{Newton-Jacobian}%
\end{equation}
where
\begin{align*}
\mathbf{A}\mathbf{=}\left[  a_{ab}\right]  ,\qquad &  a_{ab}=\int_{D} \nu\, \nabla\phi_{b}:\nabla\phi_{a}\mathbf{,}\\
\mathbf{N}^{n}=\left[  n_{ab}^{n}\right]  ,\qquad &  n_{ab}^{n}=\int%
_{D}\left(  u^{n}\cdot\nabla\phi_{b}\right)  \cdot\phi_{a},\\
\mathbf{W}^{n}=\left[  w_{ab}^{n}\right]  ,\qquad &  w_{ab}^{n}=\int%
_{D}\left(  \phi_{b}\cdot\nabla u^{n}\right)  \cdot\phi_{a}.
\end{align*}
For Picard iteration, the Newton derivative matrix~$\mathbf{W}^{n}$ is
dropped, and $\mathbf{F}^{n}=\mathbf{A}+\mathbf{N}^{n}$. The matrices are
sparse and $n_{x}$ is typically large. The divergence matrix~$\mathbf{B}$ is
defined as
\begin{equation}
\mathbf{B}=\left[  b_{cd}\right]  ,\qquad b_{cd}=\int_{D}\phi_{d}\left(
\nabla\cdot\varphi_{c}\right)  . \label{eq:B}%
\end{equation}
The residuals $\mathbf{R}^{n}$ and $\mathbf{r}^{n}$ at step~$n$ of both nonlinear iterations are given by discretization of  
(\ref{eq:residual-1})--(\ref{eq:residual-2}), and they are computed as
\begin{equation}
\left[
\begin{array}
[c]{c}%
\mathbf{R}^{n}\\
\mathbf{r}^{n}%
\end{array}
\right]  =\left[
\begin{array}
[c]{c}%
\mathbf{f}\\
\mathbf{g}%
\end{array}
\right]  -\left[
\begin{array}
[c]{cc}%
\mathbf{P}^{n} & \mathbf{B}^{T}\\
\mathbf{B} & \mathbf{0}%
\end{array}
\right]  \left[
\begin{array}
[c]{c}%
\mathbf{u}^{n}\\
\mathbf{p}^{n}%
\end{array}
\right]  , \label{eq:residual}%
\end{equation}
where $\mathbf{P}^{n}=\mathbf{A}+\mathbf{N}^{n}$ and $\mathbf{f}$ is a
discrete version of the forcing function of~(\ref{eq:NS-1}).\footnote{We use
the convention that the right-hand sides of discrete systems incorporate
Dirichlet boundary data for velocities.}

\section{Linear stability of the Navier\textendash Stokes equations}

\label{sec:stability}Following~\cite{Elman-2012-LII} let us consider, in a
general setup, the dynamical system
\begin{equation}
\mathbf{M}u_{t}=f(u,\nu), \label{eq:system-det}%
\end{equation}
where $f:%
\mathbb{R}
^{n}\times%
\mathbb{R}
\mapsto%
\mathbb{R}
^{n}$ is a nonlinear mapping, $u\in%
\mathbb{R}
^{n}$ is the state variable and $u_{t}$ is its time derivative,
$\mathbf{M}\in%
\mathbb{R}
^{n\times n}$ is the mass matrix, and $\nu$ is a parameter.
For a fixed value of~$\nu$, linear stability of the steady-state solution is
determined by the spectrum of the eigenvalue problem
\begin{equation}
\mathbf{J}v=\lambda\mathbf{M}v, \label{eq:eig-det}%
\end{equation}
where $\mathbf{J}=\frac{\partial f}{\partial u}(u(\nu),\nu)$ is the Jacobian
matrix of~$f$ evaluated at~$\nu$. The eigenvalues have a general form
$\lambda=\alpha+i\beta$, where $\alpha=\operatorname{Re}\lambda$ and
$\beta=\operatorname{Im}\lambda$,
and there are two cases: if $\alpha<0$ the perturbation decays with time, and
if $\alpha>0$ the perturbation grows.
Therefore, a change of stability can be detected by monitoring the rightmost
eigenvalues of~(\ref{eq:eig-det}).

We consider a special case of~(\ref{eq:system-det}), the time-dependent
Navier\textendash Stokes equations~(\ref{eq:NS-1})--(\ref{eq:NS-2}), 
\begin{equation}
\begin{aligned} \vec{u}_{t}&=\nu \nabla ^{2}\vec{u}-\left( \vec{u}\cdot \nabla \right) \vec{u}-\nabla p, \\ 0&=\nabla \cdot \vec{u}, \end{aligned} \label{eq:NS-time}%
\end{equation}
subject to appropriate boundary and initial conditions.
Mixed finite element discretization of~(\ref{eq:NS-time}) gives the following
Jacobian and the mass matrix, see~\cite{Elman-2012-LII}
and~\cite[Chapter~$8$]{Elman-2014-FEF} for more details,
\begin{equation}
\mathbf{J}=\left[
\begin{array}
[c]{cc}%
\mathbf{F} & \mathbf{B}^{T}\\
\mathbf{B} & \mathbf{0}%
\end{array}
\right]  \in%
\mathbb{R}
^{n_{x}\times n_{x}},\qquad\mathbf{M}=\left[
\begin{array}
[c]{cc}%
-\mathbf{G} & \mathbf{0}\\
\mathbf{0} & \mathbf{0}%
\end{array}
\right]  \in%
\mathbb{R}
^{n_{x}\times n_{x}}, \label{eq:matrices-det}%
\end{equation}
where $\mathbf{F}$ is defined as in~(\ref{Newton-Jacobian}) using the
steady-state solution of~(\ref{eq:NS-time}), $\mathbf{B}$ is defined
by~(\ref{eq:B}), and $\mathbf{G}$\ is the velocity mass matrix defined as
\[
\mathbf{G}=\left[  g_{ab}\right]  ,\qquad g_{ab}=\int_{D}\phi_{b}\,\phi_{a}.
\]
which is symmetric positive definite. Since the mass matrix$~\mathbf{M}$ is
singular, problem~(\ref{eq:eig-det}) has an infinite eigenvalue. As suggested
in~\cite{Cliffe-1994-EBM}, we replace the mass matrix$~\mathbf{M}$ with the
nonsingular, shifted mass matrix
\begin{equation}
\mathbf{M}_{\delta}=\left[
\begin{array}
[c]{cc}%
-\mathbf{G} & \delta\mathbf{B}^{T}\\
\delta\mathbf{B} & \mathbf{0}%
\end{array}
\right]  , \label{eq:M-shifted}%
\end{equation}
which maps the infinite eigenvalues of~(\ref{eq:eig-det}) to $\delta^{-1}$ and
leaves the finite ones unchanged. Then, the generalized eigenvalue
problem~(\ref{eq:eig-det}) can be replaced by
\begin{equation}
\mathbf{J}v=\lambda\mathbf{M}_{\delta}v. \label{eq:eig-det-2}%
\end{equation}
Efficient methods for estimating the rightmost pair of complex eigenvalues
of~(\ref{eq:eig-det}) (or~(\ref{eq:eig-det-2})) were studied
in~\cite{Elman-2012-LII}. Here, our goal is different. We consider
parametric uncertainty in the sense that the parameter $\nu\equiv\nu(\xi)$,
where $\xi$ is a set of random variables.

\section{The Navier\textendash Stokes equations with stochastic viscosity}

\label{sec:stochastic} Let $\left(  \Omega,\mathcal{F},\mathcal{P}\right)  $
represent a complete probability space, where $\Omega$ is the sample space,
$\mathcal{F}$ is a$~\sigma$-algebra on~$\Omega$ and $\mathcal{P}$ is a
probability measure. We will assume that the randomness in the model is
induced by a vector of independent, identically distributed (i.i.d.)\ random
variables $\xi=\left(  \xi_{1},\dots,\xi_{m_{\xi}}\right)  ^{T}$ such that
$\xi:\Omega\rightarrow\Gamma\subset\mathbb{R}^{m_{\xi}}$. Let $\mathcal{B(}%
\Gamma)$\ denote the Borel $\sigma$-algebra on$~\Gamma$ induced by$~\xi$, and
let $\rho$ 
denote the induced probability measure
for~$\xi$. The expected value of the product of measurable fuctions$~u$
and$~v$\ that depend on$~\xi$\ determines a Hilbert space~$T_{\Gamma}\equiv
L^{2}\left(  \Gamma,\mathcal{B(}\Gamma),\rho\right)  $ with inner product
\begin{equation}
\left\langle u,v\right\rangle =\mathbb{E}\left[  uv\right]  =\int_{\Gamma
}u\left(  \xi\right)  v\left(  \xi\right) \rho\,  d\xi , 
\label{eq:E}%
\end{equation}
where the symbol~$\mathbb{E}$ denotes mathematical expectation.

In computations, we 
use a finite-dimensional subspace $T_{P}\subset
T_{\Gamma}$ spanned by a set of polynomials $\left\{  \psi_{\ell}\left(
\xi\right)  \right\}  $ that are orthogonal with respect to~$\rho$, 
that is $\left\langle \psi_{k},\psi_{\ell}\right\rangle
=\delta_{k\ell}$. This is referred to as the gPC basis;
see~\cite{Ghanem-1991-SFE,Xiu-2002-WAP} for details and discussion.
For~$T_{P}$, we will use the space spanned by multivariate polynomials
in~$\{\xi_{j}\}_{j=1}^{m_{\xi}}$ of total degree~$p$, which has
dimension~$n_{\xi}={\renewcommand{\arraystretch}{0.8}\left(  \!\!%
\begin{array}
[c]{c}%
m_{\xi}+p\\
p
\end{array}
\!\!\right)  }$. We follow the setup from~\cite{Sousedik-2016-SGM} and
assume that the viscosity$~\nu$ is given by a stochastic\ expansion
\begin{equation}
\nu(\xi)=\sum_{\ell=1}^{n_{\nu}}\nu_{\ell}(x)\,\psi_{\ell}(\xi),
\label{eq:viscosity}%
\end{equation}
where $\left\{  \nu_{\ell}(x)\right\}  $\ is a set of given deterministic
spatial functions. We note that
%
this is tantamount to taking the Reynolds
number~(\ref{eq:Re}) to be stochastic.

\subsection{Stochastic linear stability and Monte Carlo simulation}

\label{sec:stochastic_stability}We are interested in a stochastic counterpart
of the generalized eigenvalue problem~(\ref{eq:eig-det-2}), that is
\begin{equation}
\mathbf{J}(\xi)v(\xi)=\lambda(\xi)\mathbf{M}_{\delta}v(\xi),
\label{eq:NS-eig-gen}%
\end{equation}
where $\mathbf{J}(\xi)$ is the nonsymmetric Jacobian matrix, which along with
the eigenvalues $\lambda(\xi)\in%
\mathbb{C}
$ and eigenvectors $v(\xi)\in%
\mathbb{C}
^{n_{x}}$ depends on the vector$~\xi$. The rightmost eigenvalue can be studied
by Monte Carlo simulation, which entails the solution
of a number of mutually independent deterministic problems
at a set of sample points$~\xi^{(i)}$, $i=1,\dots,n_{MC}$. The
sample points\ are generated randomly following the distribution of the random
variables$~\xi$, and they give realizations of the viscosity by
evaluating~(\ref{eq:viscosity}). A realization of viscosity gives rise to
deterministic functions~$\vec{u}\left(  \cdot,\xi^{\left(  i\right)  }\right)
$ and$~p\left(  \cdot,\xi^{\left(  i\right)  }\right)  $ that satisfy the
deterministic steady Navier\textendash Stokes equations, and to finite-element
approximations~$\mathbf{u}^{(i)}$ and~$\mathbf{p}^{(i)}$. The
vector~$\mathbf{u}^{(i)}$ is used to set up the Jacobian$~\mathbf{J}(\xi
^{(i)})$ and solving~(\ref{eq:NS-eig-gen}) provides a realization of the
rightmost eigenvalue$~\lambda(\xi^{(i)})$.

In Monte Carlo simulation this procedure is thus performed for every sample
$i=1,\dots,n_{MC}$, and the moments of the eigenvalue are obtained from
ensemble averaging. We will also use the term \emph{simulator} and denote it
by$~\eta$ for the computer code computing the rightmost eigenvalue
of~(\ref{eq:NS-eig-gen}) for given input parameters$~\xi$. Since use of the
simulator is in general computationally expensive,
we are interested in construction of an \emph{emulator}, which is a
computationally cheap \emph{surrogate} of the full model that can be easily
evaluated for any value of the input parameters. We will denote use of an
emulator by $\lambda_{\star}(\xi)=\eta_{\star}(\xi)$, where the symbol$~\star
$\ stands for any of the three approaches to emulation and surrogate
construction discussed next.


\subsection{Polynomial chaos surrogate}

\label{sec:collocation}

Both Monte Carlo and stochastic collocation methods are based on sampling. For
stochastic collocation, the sample points$~\xi^{(q)}$, $q=1,\dots,n_{q}$,
consist of a set of predetermined \emph{collocation points}. This approach
derives from a methodology for performing quadrature or interpolation in
multidimensional space using a small number of points, a so-called sparse
grid~\cite{Gerstner-1998-NIU,Novak-1996-HDI}. There are two ways to
implement stochastic collocation, either by constructing a Lagrange
interpolating polynomial, or, in the so-called pseudospectral approach, by
performing a discrete projection into~$T_{P}$%
~\cite{Babuska-2010-SCM,Xiu-2010-NMS}. We use the second approach. In
particular, we will search for expansions of the eigenvalue ${\lambda(\xi)}$
in the form
\begin{equation}
\lambda(\xi)=\sum_{k=1}^{n_{\xi}}\lambda_{k}\psi_{k}(\xi), \label{eq:sol_mat}%
\end{equation}
where $\lambda_{k}\in\mathbb{%
\mathbb{C}
}$ are coefficients corresponding to the basis$~\left\{  \psi_{k}\right\}  $
defined by a discrete projection
\begin{equation}
\lambda_{k}=\left\langle \lambda,\psi_{k}\right\rangle ,\qquad k=1,\dots
,n_{\xi}. \label{eq:sol_mat-proj}%
\end{equation}
The coefficients in~(\ref{eq:sol_mat}) are determined by
evaluating~(\ref{eq:sol_mat-proj}) (see~(\ref{eq:E})), using numerical
quadrature as
\begin{equation}
\lambda_{k}=\sum_{q=1}^{n_{q}}\lambda(\mathbb{\xi}^{(q)})\psi_{k}(\mathbb{\xi
}^{(q)})w^{(q)}, \label{eq:quadrature}%
\end{equation}
where $\mathbb{\xi}^{(q)}$ are the quadrature (collocation) points and
$w^{(q)}$ are quadrature weights. That is, the evaluations of coefficients
in~(\ref{eq:sol_mat-proj}) entail solving a set of independent deterministic
eigenvalue problems at a set of sample points. Details of the rule we use in
our numerical experiments are discussed in Section~\ref{sec:numerical}, and we
refer, e.g., to monograph~\cite{LeMaitre-2010-SMU} for more details.

Once the coefficients in~(\ref{eq:sol_mat-proj}) have been determined, the
stochastic collocation emulator~$\eta_{\operatorname*{SC}}$ is
\begin{equation}
\lambda_{\operatorname*{SC}}(\xi)=\eta_{\operatorname*{SC}}\left(  \xi\right)
=\sum_{k=1}^{n_{\xi}}\lambda_{k}\psi_{k}(\xi). \label{eq:gPC-surrogate}%
\end{equation}
See~\cite{Andreev-2012-STA} for analysis showing convergence of this approximation 
for self-adjoint problems.

\subsection{Surrogate based on Gaussian process regression}

\label{sec:GPR} In Gaussian process regression we assume that if the process
depends on~$n_{\star}$ inputs in$~m_{\xi}$ dimensions, then the output is
an~$n_{\star}$-dimensional vector. Specifically, the output is modeled as
\begin{equation}
\lambda_{\operatorname*{GP}}\left(  \xi\right)  =\eta_{\operatorname*{GP}%
}\left(  \xi\right)  =\mu+z(\xi), \label{eq:GP}%
\end{equation}
where we consider $\mu$ as a constant, which is also common in practice, 
and~$z$ is a Gaussian process to be determined. The distribution of the output is
multivariate normal with mean~$\mu$. For the covariance function$~R$ we
consider the so-called squared exponential kernel function, and we note that
it is proportional to a correlation (or kernel) matrix$~C$ by a constant of
proportionality$~\sigma_{f}^{2}$ called the variance ($\sigma_{f}$ is the
standard deviation) via $R=\sigma_{f}^{2}C$. Specifically, the correlation
function$~C$ has the entries given by
\[
C\left(  \xi,\xi^{\prime}\right)  =\exp\left[  -\frac{1}{2}\frac{\left(
\xi-\xi^{\prime}\right)  ^{T}\left(  \xi-\xi^{\prime}\right)  }{\sigma_{\ell}%
}\right]  ,
\]
where $\sigma_{\ell}$ is the correlation length.
The prior for the simulator
is
\[
\eta_{\operatorname*{GP}}^{\text{prior}}\left(  \xi\right)  \sim
\mathcal{N}\left(  \mu,R\left(  \xi,\xi\right)  \right)  ,
\]
where $\mathcal{N}$ denotes the multivariate normal distribution. The
parameters $\mu$, $\sigma_{f}$, and$~\sigma_{\ell}$\ are estimated from the
simulator runs at the experimental design points $\xi^{(t)}$, $t=1,\dots
,n_{d}$, with results collected in a vector$~\lambda_{GP}^{d}$. Let us define
the correlation matrix $C_{d}$ with entries $c_{ij}=C(\xi_{i},\xi_{j})$, where
$i,j=1,\dots,n_{d}$, and let us denote by $H$ a vector of ones with
length$~n_{d}$. 
Assuming a standard noninformative prior for variance parameters following~\cite{Owen-2017-CSB}, 
we estimate
\begin{align*}
\hat{\mu}  &  =\left(  H^{T}C_{d}^{-1}H\right)  ^{-1}H^{T}C_{d}^{-1}\lambda
_{\operatorname*{GP}}^{d},\\
\hat{\sigma}_{f}  &  =\left(  \lambda_{\operatorname*{GP}}^{d}-\hat{\mu
}H\right)  ^{T}C_{d}^{-1}\left(  \lambda_{\operatorname*{GP}}^{d}-\hat{\mu
}H\right)  .
\end{align*}
The correlation length is estimated by maximizing the logarithm of the
likelihood$~L$~as
\[
\hat{\sigma}_{\ell}=\arg\max_{\sigma_{\ell}}\left[  \log L\left(  \sigma
_{\ell}|\lambda_{\operatorname*{GP}}^{d}\right)  \right]  ,
\]
where the likelihood for the correlation length is
\[
L\left(  \sigma_{\ell}|\lambda_{\operatorname*{GP}}^{d}\right)  \propto\left(
\hat{\sigma}_{f}^{2}\right)  ^{-\left(  n_{d}-n_{\mu}\right)  /2}\left\vert
C_{d}\right\vert ^{-1/2}\left\vert H^{T}C_{d}^{-1}H\right\vert ^{-1/2},
\]
where $\left\vert \cdot\right\vert$ is the determinant, and we use
$n_{\mu}=1$ since we consider constant $\mu$ in~(\ref{eq:GP}).

After the parameters have been determined, the Gaussian process emulator~$\eta
_{\operatorname*{GP}}$ is specified by a posterior distribution which is a
Student's $t$-distribution with $n_{d}-n_{\mu}$ degrees of freedom
\begin{equation}
\eta_{\operatorname*{GP}}\left(  \xi\right)  \sim t_{n_{d}-n_{\mu}}\left(
M^{\ast}(\xi)|R^{\ast}(\xi,\xi)\right)  . \label{eq:GP-post}%
\end{equation}
The posterior mean and covariance functions in (\ref{eq:GP-post}) are defined,
respectively, as
\begin{align*}
M^{\ast}(\xi)  &  =\hat{\mu}+\hat{R}(\xi)C_{d}^{-1}\left(  \lambda
_{\operatorname*{GP}}^{d}-\hat{\mu}H\right) \\
R^{\ast}(\xi,\xi^{\prime})  &  =\frac{\hat{\sigma}_{f}^{2}}{n_{d}-n_{\mu}%
-2}\left[  C(\xi,\xi^{\prime})-\hat{R}(\xi)C_{d}^{-1}\hat{R}(\xi^{\prime}%
)^{T}+Q(\xi)\left(  H^{T}C_{d}^{-1}H\right)  ^{-1}Q(\xi^{\prime})^{T}\right]
,
\end{align*}
where $\hat{R}(\xi)$ is a (row) vector of correlations between $\xi$ and the
experimental design points, and $Q(\xi)=1-\hat{R}(\xi)C_{d}^{-1}H$. In
implementation, we use the \textsc{Matlab} functions \texttt{fitrgp}, and
\texttt{predict} with more details given in discussion of numerical
experiments in Section~\ref{sec:numerical}. We also note that even though the
emulator$~\eta_{\operatorname*{GP}}$ readily provides uncertainty information
through the posterior distribution~(\ref{eq:GP-post}),
we explore $~\eta_{\operatorname*{GP}}$ by evaluating it directly so that it
is treated in a manner consistent with the other emulators $~\eta
_{\operatorname*{SC}}$ and $~\eta_{\operatorname*{NN}}$, the latter of which
is discussed next.

\subsection{Neural network surrogate}

The final surrogate is based on a shallow (as opposed to deep) neural network
with a single hidden layer and hyperbolic tangent sigmoid transfer function
\texttt{tansig}, which is mathematically equivalent to $\tanh$, 
see~\cite{Vogl-1988-ACB}. The goal is to develop an emulator
\[
\lambda_{\operatorname*{NN}}\left(  \xi\right)  =\eta_{\operatorname*{NN}%
}\left(  \xi\right)  ,
\]
based on nonlinear regression and supervised learning. The network is trained
as follows. We are given a training set of inputs and targets in the form
$\left\{  \xi^{(t)},\lambda(\xi^{(t)})\right\}  $, $t=1,\dots,n_{t}$, and the
training data is split into groups used for training, testing and validation.
The neural network emulator$~\eta_{\operatorname*{NN}}$ is initialized
randomly, and the task of the training is to produce a network that produces
small errors on the training set but also responds well to additional inputs.
In that case we say that the network generalizes well. The process of training
a neural network entails tuning the values of the weights and biases of the
network to optimize network performance by minimizing the sum of squared
errors
\[
\frac{1}{n_{t}}\sum_{t=1}^{n_{t}}\left(  \lambda(\xi^{(t)})-\lambda
_{\operatorname*{NN}}(\xi^{(t)})\right)  ^{2}.
\]
The specific algorithm we use for the training is the Bayesian regularization
backpropagation, in which the weight and bias values are updated according to
Levenberg\textendash Marquardt optimization, 
see~\cite{DanForesee-1997-GNA,MacKay-1992-BI} for details. In implementation, we
use \textsc{Matlab} functions \texttt{fitnet}, \texttt{train} and \texttt{net}
with more details given in Section~\ref{sec:numerical}.

\subsection{Validation and assessment of the surrogate models}

\label{sec:validation} After the surrogates are built, we would like to assess
and compare their quality. Our strategy is similar to that used 
by~\cite{Owen-2017-CSB}. Specifically, for the validation of the surrogates
constructed using the emulators we used Monte Carlo simulation, for which the
input parameters $\xi^{(i)}$, $i=1,\dots,n_{MC}$, are distinct from the input
parameters used to build the surrogates. The validation metric is then given
by root mean square error ($\operatorname*{RMSE}$) defined as
\[
\operatorname*{RMSE}=\sqrt{\frac{1}{n_{MC}}\sum_{i=1}^{n_{MC}}\left(
\lambda_{\star}(\xi^{(i)})-\lambda(\xi^{(i)})\right)  ^{2}},
\]
where the symbol $\star$ denotes any of the SC, GP or NN emulators. We used
the Monte Carlo sample points$~\xi^{(i)}$, $i=1,\dots,n_{MC}$. Since
$\operatorname*{RMSE}$ represents the distance between a surrogate and the
Monte Carlo simulator across the input prameters space, low
$\operatorname*{RMSE}$ values are favorable.

Next, we compute the mean and variance of each surrogate, $\mu_{\star}$ and
$\sigma_{\star}$, respectively, and we estimate those provided by the
emulators using empirical formulas given as
$$
\mu_{\star}   =\frac{1}{n_{MC}}\sum_{i=1}^{n_{MC}}\lambda_{\star}(\xi^{(i)}), \qquad
\sigma_{\star}   =\sqrt{\frac{1}{n_{MC}}\sum_{i=1}^{n_{MC}}\left(
\lambda_{\star}(\xi^{(i)})-\mu_{\star})\right)  ^{2}}.
$$
Although for stochastic collocation both quantities above could be calculated
directly from the gPC\ coefficients, here we used the above formulas also
with$~\eta_{\operatorname*{SC}}$. Since we want to detect 
instability, we also use the surrogates to estimate the probability that the
rightmost eigenvalue is nonnegative as
\[
\Pr\left(  \lambda_{\star}\geq0\right)  \approx\frac{1}{n_{MC}}\sum
_{i=1}^{n_{MC}}\mathbf{1}\left(  \lambda_{\star}(\xi^{(i)})\geq0\right)  ,
\]
where $\mathbf{1}$ denotes the indicator (1 or $0$) function. Finally, we also
test the ability of the surrogate to reconstruct the probability density
function of the simulator output, which we do using a kernel density estimator
with Gaussian kernel provided by the \textsc{Matlab}\ function
\texttt{ksdensity}.

\begin{remark}We note that only the one of these, the neural network emulator, exactly fits within the paradigm of ``machine learning'' methods in the sense that it constructs a neural network.  However, we view all of them as methods based on learning, in the sense that the surrogate is built from data obtained from a training set,  where for stochastic collocation the learning process is the construction of the solution at the collocation points, and for Gaussian process regression, it is the construction of the mean, variance and correlation length from the simulation at the design points.
\end{remark}

\section{Numerical experiments}

\label{sec:numerical}We implemented the Navier\textendash Stokes solver
in~\textsc{Matlab} version 9.7.0.1190202 (R2019b)\ using the
\textsc{IFISS~3.5} package~\cite{ERS-SIREV}, and we tested the simulator and
the emulators using two benchmark problems: flow around an obstacle and an
expansion flow around a symmetric step. 
These are representative examples that exhibit important types of 
bifurcation, a Hopf bifurcation for the first (where the critical eigenvalues are a complex
conjugate pair) 
and a pitchfork bifurcation for the second (with a real critical eigenvalue)~\cite{Cliffe-Spence-Taverner,Govaerts}.
For both examples, we consider perturbations of mean viscosities that are near values
leading to bifurcations.

\begin{figure}[ptbh]
\begin{center}
\includegraphics[width=8.8cm]{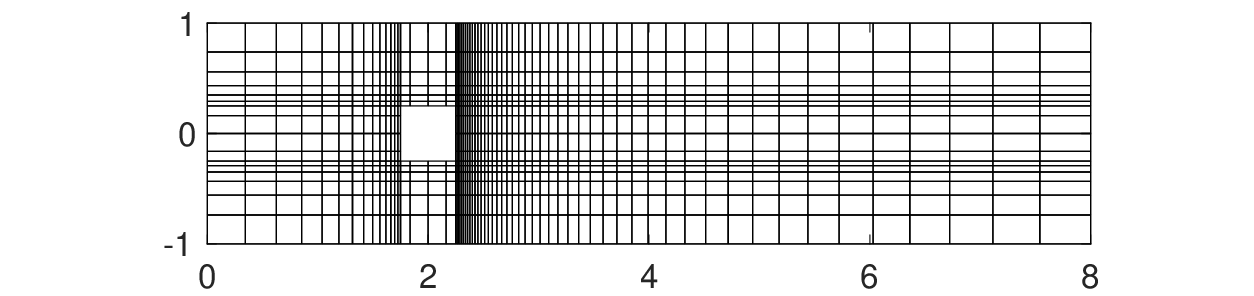}
\end{center}
\vspace{-.2in}
\caption{Finite element mesh for the flow around an obstacle problem. }%
\label{fig:mesh-obstacle}%
\end{figure}

\begin{figure}[ptbh]
\centering
\begin{tabular}
[c]{cc}%
\includegraphics[width=6.2cm]{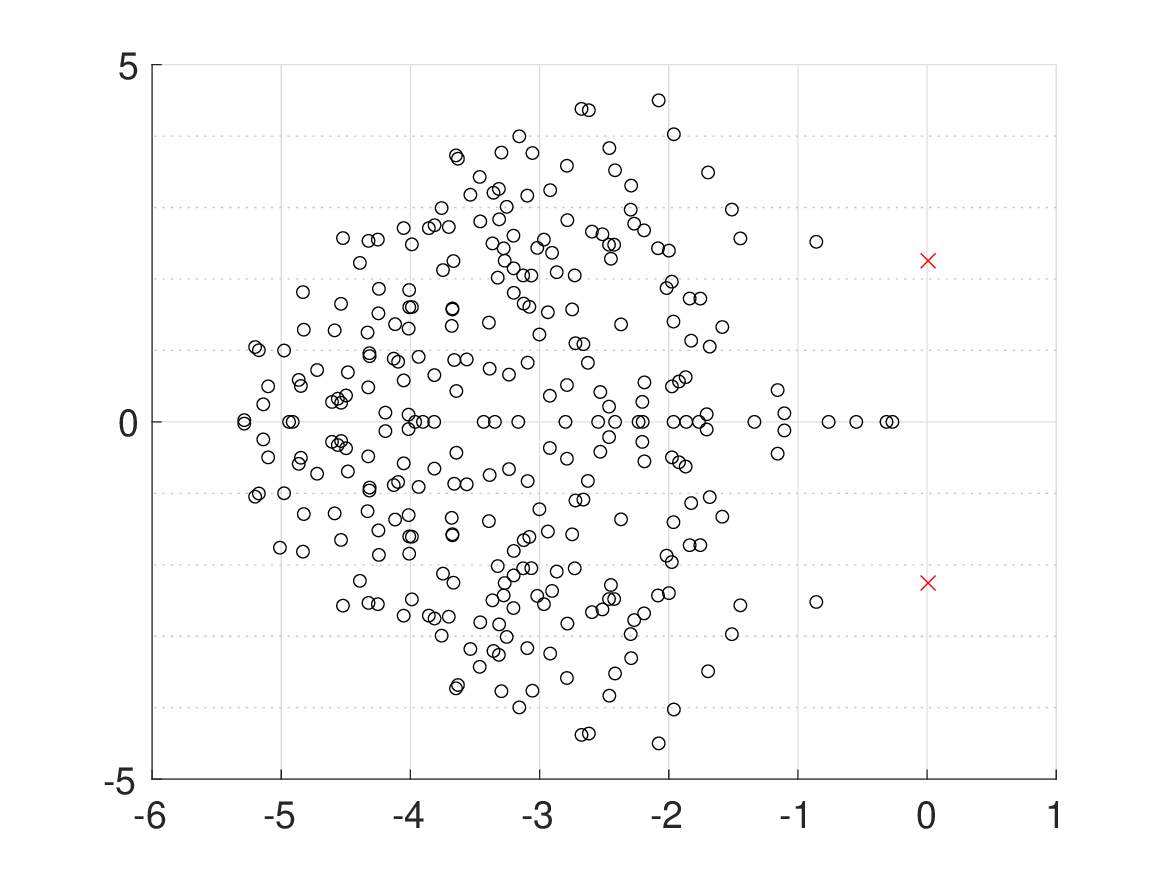} &
\includegraphics[width=6.2cm]{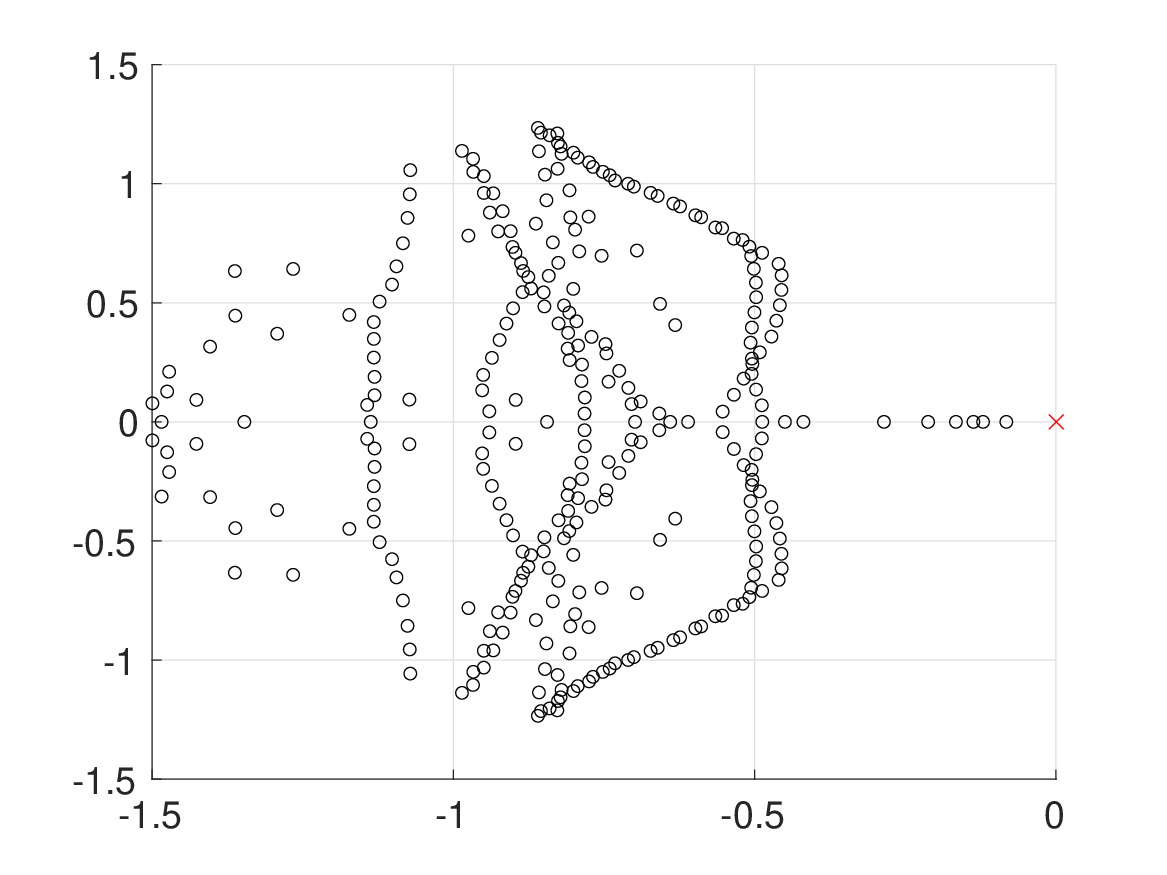}
\end{tabular}
\caption{An image of the complex plane and $300$ eigenvalues with the largest real part of the deterministic
eigenvalue problem with mean viscosity (i.e., $\nu=\nu_{1}$
in~(\ref{eq:viscosity})) for the two examples: flow around an obstacle (left)
and expansion flow around a symmetric step (right). The rightmost eigenvalues
are indicated by a red cross. }%
\label{fig:lambdaM}%
\end{figure}

For the solution of the steady
Navier\textendash Stokes problem in the simulator we used a hybrid strategy in which an
initial approximation is obtained from solution of the stochastic Stokes
problem, after which several steps of Picard iteration are used to improve the
solution, followed by Newton iteration.
The convergence test was for
the Euclidean norm of the algebraic residual~(\ref{eq:residual}) to satisfy 
\[
\left\Vert \left[
\begin{array}
[c]{c}%
\mathbf{R}^{n}\\
\mathbf{r}^{n}%
\end{array}
\right]  \right\Vert \leq10^{-8}\left\Vert \left[
\begin{array}
[c]{c}%
\mathbf{f}\\
\mathbf{g}%
\end{array}
\right]  \right\Vert .
\]
Next, the eigenvalue problems~(\ref{eq:eig-det-2}), in which $\mathbf{M}%
_{\delta}$ is defined by~(\ref{eq:M-shifted}) with $\delta=-10^{-2}$ as
in~\cite{Elman-2012-LII}, were solved using the function \texttt{eigs} in
\textsc{Matlab}. The $300$ eigenvalues with the largest real part of the
deterministic eigenvalue problem with mean viscosity~$\nu_{1}$ for each of the
two examples are displayed in Figure~\ref{fig:lambdaM}. The
viscosity~(\ref{eq:viscosity}) is parameterized using $m_{\xi}=2$ random
variables. For the Monte Carlo method we used $10^{3}$ sample points generated
randomly following the
distribution of the random variables $\xi$. For stochastic collocation\ we
used Smolyak sparse grid
and grid level$~4$. With these settings, there were $n_{q}=29$ points on the
sparse grid, and this set of quadrature points was used to design all three
emulators $\eta_{\operatorname*{SC}}$, $\eta_{\operatorname*{GP}}$\ and
$\eta_{\operatorname*{NN}}$, that is $n_{q}=n_{d}=n_{t}$. For the
GP\ regression (and also for the training of the neural network) we
standardize the data before the regression. To this end let $\mu^{d}$ and
$\sigma^{d}$\ denote the mean and standard deviation of the rightmost
eigenvalues$~\lambda(\xi^{(q)})$ calculated using the simulator at the
quadrature points~$\xi^{(q)}$, $q=1,\dots,n_{q}$. The data points passed to
the GP\ regression function \texttt{fitrgp} in \textsc{MATLAB} are scaled as
\begin{equation}
\lambda(\xi^{(q)})\leftarrow\frac{\lambda(\xi^{(q)})-\mu^{d}}{\sigma^{d}%
},\qquad q=1,\dots,n_{q}, \label{eq:training-scale}%
\end{equation}
and the results$~\lambda_{\operatorname*{GP}}(\xi)$ of the emulator function
\texttt{predict} are descaled as
\begin{equation}
\lambda_{\operatorname*{GP}}(\xi)\leftarrow\sigma^{d}\,\lambda
_{\operatorname*{GP}}(\xi)+\mu^{d}. \label{eq:training-descale}%
\end{equation}
For the neural network emulator we use function \texttt{fitnet} in
\textsc{MATLAB} to construct a neural network with one hidden layer of
$20$\ neurons, and we set the training algorithm to use Bayesian
regularization. The training parameters used in the actual training function
\texttt{train} are divided to $80\%$ for training, $10\%$ for testing and
$10\%$ for validation. While we do not have a general strategy to find the
optimal size of the neural network, we empirically\ tried to find as small
a\ network as possible that would still match the Monte Carlo simulation
reasonably well. We used scaling~(\ref{eq:training-scale}) for the training,
and descaling~(\ref{eq:training-descale}) for the emulator predictions given
by the function \texttt{net} in \textsc{MATLAB}.



\begin{figure}[ptbh]
\begin{center}
\includegraphics[width=13cm]{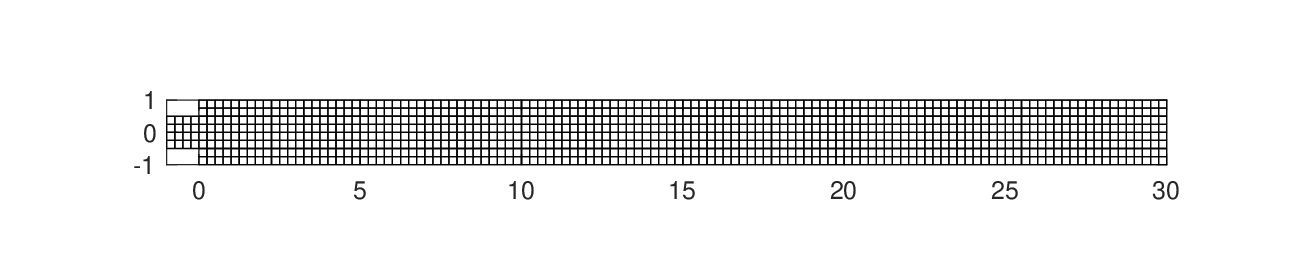}
\end{center}
\vspace{-.4in}
\caption{Finite element mesh for the expansion flow around a symmetric step. }%
\label{fig:mesh-symstep}%
\end{figure}

\subsection{Flow around an obstacle}
\label{sec:obstacle}

\begin{figure}[ptbh]
\centering
\begin{tabular}
[c]{cc}%
\includegraphics[width=6.2cm]{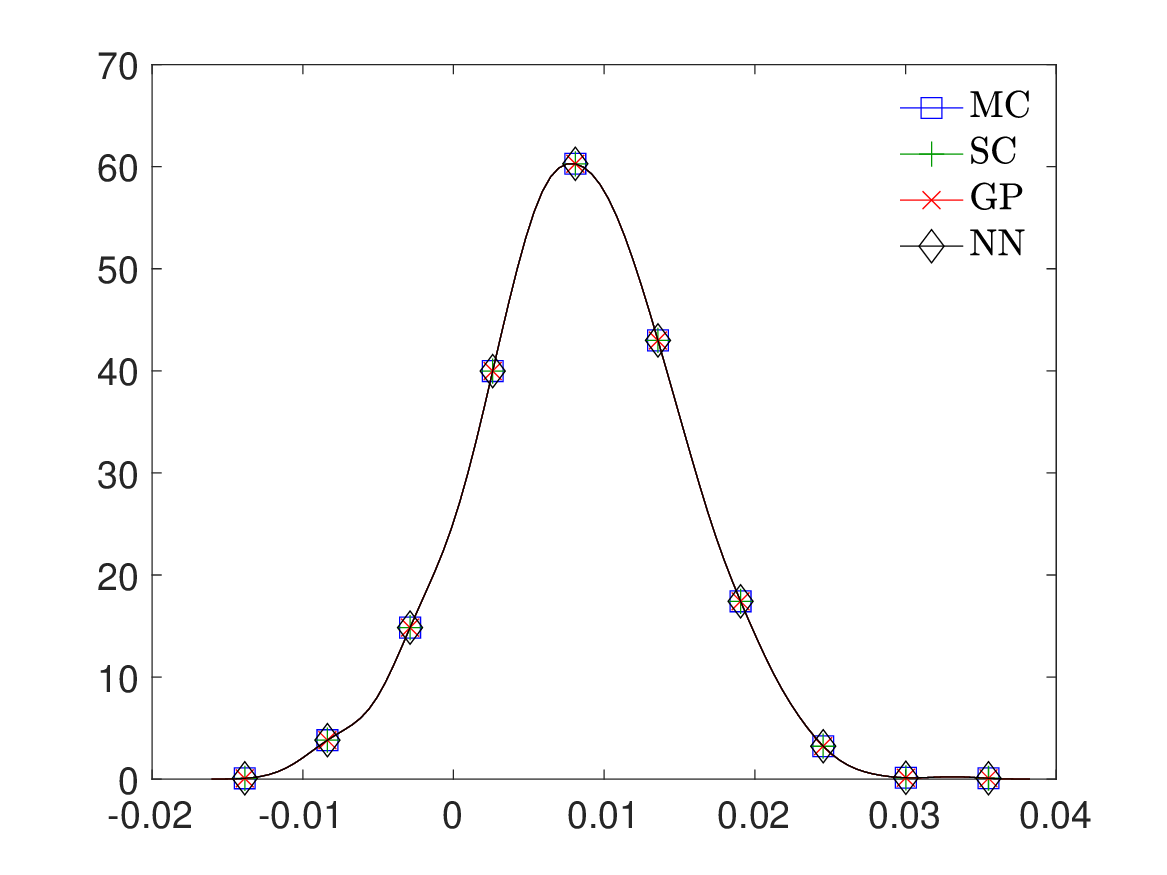} &
\includegraphics[width=6.2cm]{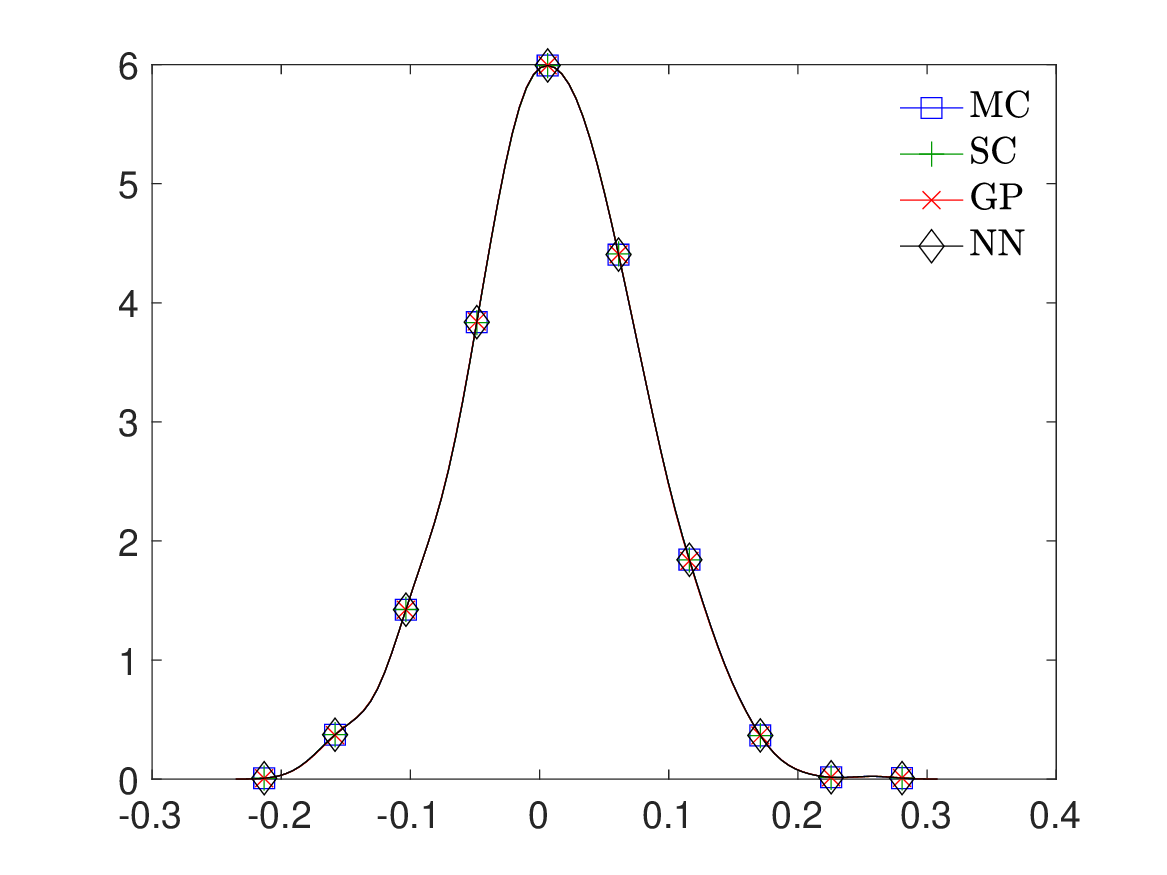}
\end{tabular}
\caption{Plots of the pdf estimate of the real part of the rightmost
eigenvalue obtained using Monte Carlo (MC), stochastic collocation (SC),
Gaussian process regression (GP) and neural network (NN) for the flow around
an obstacle with $CoV=1\%$ (left) and $CoV=10\%$ (right).}%
\label{fig:obstacle-pdf-Re}%
\end{figure}

For the first example, we consider flow around an obstacle in a similar setup
as studied in~\cite{Sousedik-2016-SGM}.
The domain of the channel and the discretization are shown in
Figure~\ref{fig:mesh-obstacle}. The spatial discretization uses a stretched
grid with $1008$ \textbf{\textit{Q}}$_{2}-$\textbf{\textit{Q}}$_{1}$
(\emph{Taylor\textendash Hood})
finite\ elements. 
There are $8416$ velocity and $1096$ pressure degrees of freedom. The
viscosity$~\nu(x,\xi)$\ was taken to be a truncated lognormal process
transformed from an underlying Gaussian process~\cite{Ghanem-1999-NGS}.
That is, $\psi_{\ell}(\xi)$, $\ell=1,\dots,n_{\nu}$, is a set of Hermite
polynomials, which also specifies the expansion of
viscosity~(\ref{eq:viscosity}) used in the simulator. Denoting the
coefficients of the Karhunen\textendash Lo\`{e}ve expansion of the Gaussian process
by$~g_{j}(x)$ and $\iota_{j}=\xi_{j}-g_{j},$ $j=1,\dots,m_{\xi}$, the
coefficients in expansion~(\ref{eq:viscosity}) are computed as
\[
\nu_{\ell}(x)=\frac{\mathbb{E}\left[  \psi_{\ell}(\iota)\right]  }%
{\mathbb{E}\left[  \psi_{\ell}^{2}(\iota)\right]  }\exp\left[  g_{0}+\frac
{1}{2}\sum_{j=1}^{m_{\xi}}\left(  g_{j}(x)\right)  ^{2}\right]  .
\]
The covariance function of the Gaussian process, for points $X_{1}=(x_{1}%
,y_{1})$ and $X_{2}=(x_{2},y_{2})$ in$~D$, was chosen to be
\begin{equation}
C_{\text{rf}}\left(  X_{1},X_{2}\right)  =\sigma_{g}^{2}\exp\left(
-\frac{\left\vert x_{2}-x_{1}\right\vert }{L_{x}}-\frac{\left\vert y_{2}%
-y_{1}\right\vert }{L_{y}}\right)  , \label{eq:covariance_kernel}%
\end{equation}
where$~L_{x}$ and $L_{y}$ are the correlation lengths of the random variables
$\xi_{i}$, $i=1,\dots,m_{\xi}$, in the $x$ and $y$\ directions, respectively,
and $\sigma_{g}$ is the standard deviation of the Gaussian random field. The
correlation lengths were set to be equal to $25\%$ of the width and height of
the domain. The coefficient of variation $CoV$ of the lognormal field, defined
as $CoV=\sigma_{\nu}/\nu_{1}$, where $\sigma_{\nu}$ is the standard deviation
and $\nu_{1}$ is the mean viscosity, was $1\%$ or $10\%$. According
to~\cite{Matthies-2005-GML}, in order to guarantee a complete representation
of the lognormal process by~(\ref{eq:viscosity}) the degree of polynomial
expansion of~$\nu(x,\xi)$ should be twice the degree of the expansion of the
solution. We follow the same strategy here. Therefore, the values of $n_{\xi}$
and $n_{\nu}$ are, see, e.g.~\cite[p.~87]{Ghanem-1991-SFE}
or~\cite[Section~5.2]{Xiu-2010-NMS}, $n_{\xi}=\frac{\left(  m_{\xi}+p\right)
!}{m_{\xi}!p!}$, $n_{\nu}=\frac{\left(  m_{\xi}+2p\right)  !}{m_{\xi}!\left(
2p\right)  !}$. For the gPC expansion of eigenvalues~(\ref{eq:sol_mat}), the
maximal degree of gPC\ expansion is $p=3$, so then $n_{\xi}=10$ and $n_{\nu
}=28$. We assumed that the random variables $\left\{  \xi_{\ell}\right\}
_{\ell=1}^{m_{\xi}}$ follow a normal distribution and used Smolyak sparse grid
with Gauss\textendash Hermite quadrature points for collocation. For the solution of the
Navier\textendash Stokes problem we used the hybrid strategy with $6$ steps of Picard
iteration followed by at most $15$ steps of Newton iteration. We used mean
viscosity $\nu_{1}=5.36193\times10^{-3}$, which corresponds to Reynolds number
$Re=373$, and the rightmost eigenvalue pair is $0.0085\pm2.2551i$, see the
left panel in Figure~\ref{fig:lambdaM}. Table~\ref{tab:obstacle} presents the
results of validation and assessment of the surrogates using the indicators
from Section~\ref{sec:validation}. It is evident that for both $CoV$
$1\%$\ and $10\%$ the values of $\operatorname*{RMSE}$ are small for all
surrogates with the smallest value for the stochastic collocation, where we
note that we used the same values of~$\xi^{(i)}$ in the Monte Carlo simulation
and also for sampling the gPC\ surrogate~(\ref{eq:gPC-surrogate}). All values
of $\mu$ and $\sigma$ are in close agreement, and in particular, all values of
$\operatorname*{RMSE}$ are smaller than the corresponding values of $\mu$ (and
$\sigma$) by at least two orders of magnitude. Also, all emulators indicate
reliably the probability of the rightmost eigenvalue being nonnegative.
Finally, Figure~\ref{fig:obstacle-pdf-Re} displays the probability density
function (pdf) estimates of the rightmost eigenvalue.
The estimates were obtained using \textsc{Matlab} function \texttt{ksdensity} 
for sampled gPC\ expansions. In all cases, we see an excellent agreement of
the plots in the left panel corresponding to $CoV=1\%$ and in the right panel
corresponding to $CoV=10\%$.\ 

\begin{table}[b]
\caption{Flow around an obstacle: validation of the surrogate models by Monte
Carlo (MC) simulation using root mean square error (RMSE), their assessment
using estimates of the mean~$\mu$, standard deviation~$\sigma$, and
probability that the rightmost eigenvalue is nonnegative. The surrogates are
based on stochastic collocation (SC), Gaussian process regression (GP) and
neural network (NN), and the measures are defined in
Section~\ref{sec:validation}.}%
\label{tab:obstacle}
\begin{center}
\begin{tabular}
[c]{|c|c|c|c|c|}\hline
& MC & SC & GP & NN\\\hline
\multicolumn{5}{|c|}{$CoV=1\%$}\\\hline
$\operatorname{RMSE}$ & - & $4.1859\times10^{-8}$ & $2.1709\times10^{-6}$ &
$5.0301\times10^{-7}$\\\hline
$\mu$ & $8.3579\times10^{-3}$ & $8.3579\times10^{-3}$ & $8.3571\times10^{-3}$
& $8.3579\times10^{-3}$\\
$\sigma$ & $6.5356\times10^{-3}$ & $6.5356\times10^{-3}$ & $6.5355\times
10^{-3}$ & $6.5356\times10^{-3}$\\\hline
$\operatorname{Pr}(\lambda\geq0)$ & \multicolumn{4}{|c|}{$89.8\%$}\\\hline
\multicolumn{5}{|c|}{$CoV=10\%$}\\\hline
$\operatorname{RMSE}$ & - & $9.4232\times10^{-5}$ & $3.9827\times10^{-4}$ &
$2.9063\times10^{-5}$\\\hline
$\mu$ & $1.1279\times10^{-2}$ & $1.1277\times10^{-2}$ & $1.1235\times10^{-2}$
& $1.1277\times10^{-2}$\\
$\sigma$ & $6.5819\times10^{-2}$ & $6.5818\times10^{-2}$ & $6.5789\times
10^{-2}$ & $6.5813\times10^{-2}$\\\hline
$\operatorname{Pr}(\lambda\geq0)$ & \multicolumn{2}{|c|}{$56.5\%$} &
\multicolumn{1}{|c|}{$56.4\%$} & \multicolumn{1}{|c|}{$56.5\%$}\\\hline
\end{tabular}
\end{center}
\end{table}

\subsection{Expansion flow around a symmetric step}

\begin{figure}[ptbh]
\centering
\begin{tabular}
[c]{cc}%
\includegraphics[width=6.2cm]{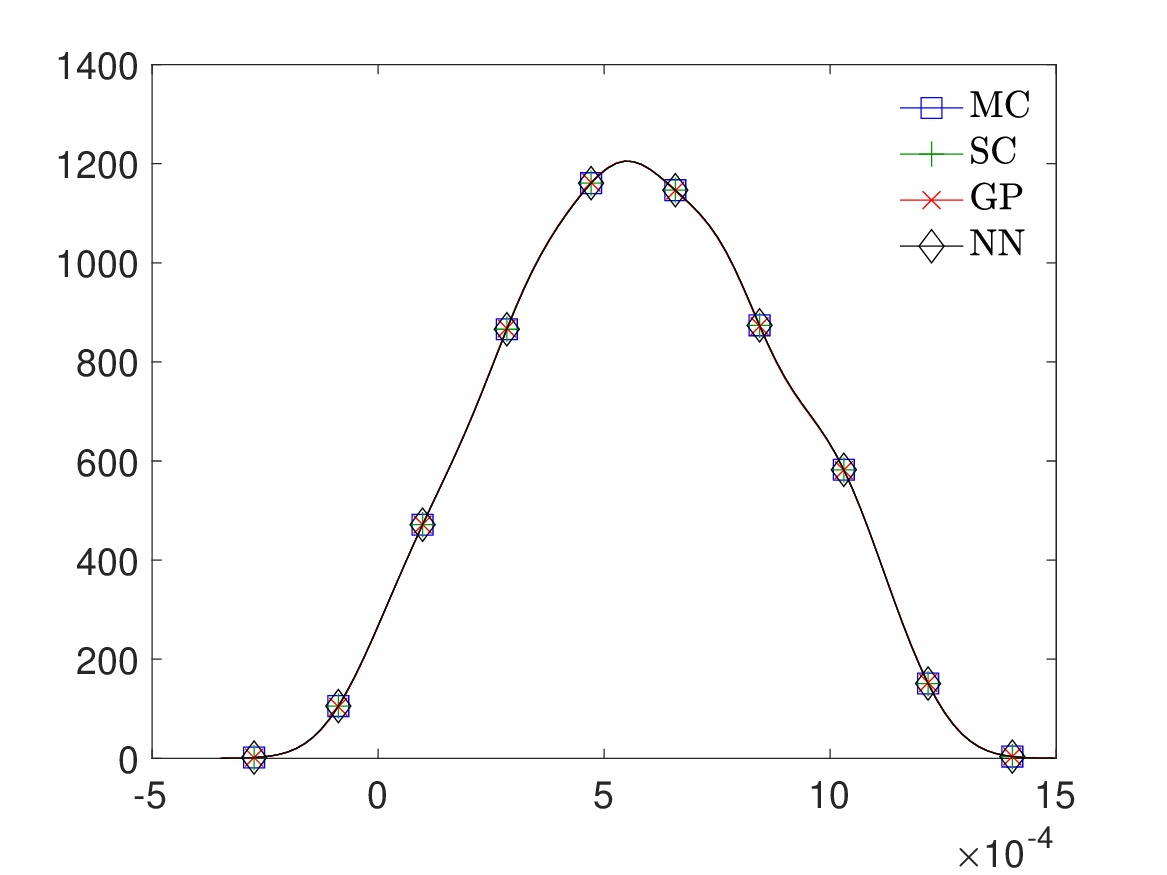} &
\includegraphics[width=6.2cm]{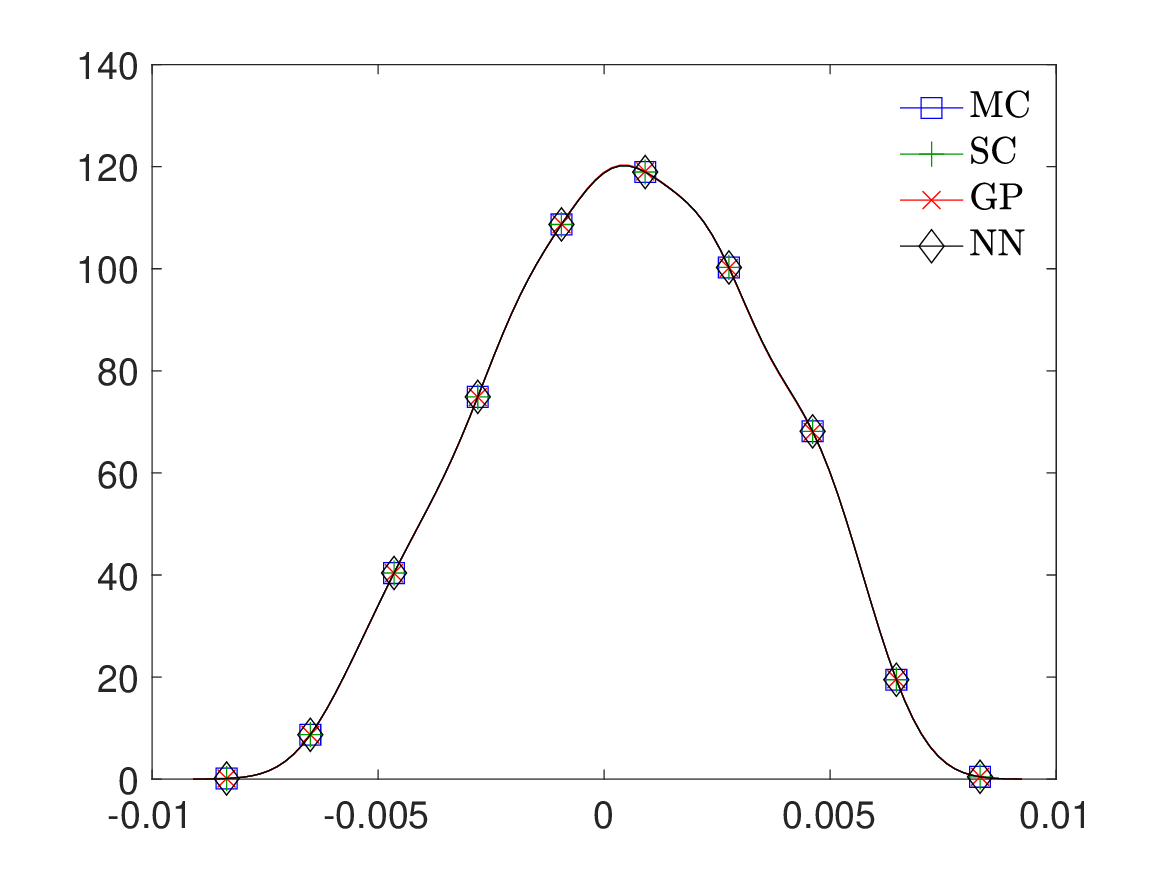}
\end{tabular}
\caption{Plots of the pdf estimate of the real part of the rightmost
eigenvalue obtained using Monte Carlo (MC), stochastic collocation (SC),
Gaussian process regression (GP) and neural network (NN) for the expansion
flow around a symmetric step with $CoV=1\%$ (left) and $CoV=10\%$ (right).}%
\label{fig:symstep-pdf-Re}%
\end{figure}

For the second example, we consider an expansion flow around a symmetric step.
The domain and its discretization are shown in Figure~\ref{fig:mesh-symstep}.
The spatial discretization uses a uniform grid with $976$ \textbf{\textit{Q}%
}$_{2}-$\textbf{\textit{P}}$_{-1}$ finite\ elements, which provide a stable
discretization for the rectangular grid \cite[p.~139]{Elman-2014-FEF}.
There are $8338$ velocity and $2928$ pressure degrees of freedom. For the
viscosity we considered a random field with affine dependence on the random
variables $\xi$ given as
\begin{equation}
\nu(x,\xi)=\nu_{1}+\sigma_{\nu} \textstyle{\sum_{\ell=2}^{n_{\nu}}} \nu_{\ell}(x)\,\xi
_{\ell-1}, \label{eq:viscosity-KL}%
\end{equation}
where $\nu_{1}$ is the mean and $\sigma_{\nu}=CoV\cdot\nu_{1}$ the standard
deviation of the viscosity, $n_{\nu}=m_{\xi}+1$, and $\nu_{\ell+1}%
=\sqrt{3\lambda_{\ell}}v_{\ell}(x)$ with $\left\{  \left(  \lambda_{\ell
},v_{\ell}(x)\right)  \right\}  _{\ell=1}^{m_{\xi}}$ are the eigenpairs of the
eigenvalue problem associated with the covariance kernel of the random field.
As in the previous example, we used the values $CoV=1\%$, and $10\%$. We
considered the covariance kernel~(\ref{eq:covariance_kernel}),
with correlation lengths set to $12.5\%$ of the width and $25\%$\ of
the height of the domain. We assumed that the random variables $\left\{
\xi_{\ell}\right\}  _{\ell=1}^{m_{\xi}}$ follow a uniform distribution over
$(-1,1)$. 
Note that~(\ref{eq:viscosity-KL}) can be viewed as a special case
of~(\ref{eq:viscosity}), which consists of only linear terms of$~\xi$. For the
parametrization of viscosity by~(\ref{eq:viscosity-KL}), which then specifies
the simulator, we used the same stochastic dimension~$m_{\xi}$ and degree of
polynomial expansion$~p$ as in the previous example: $m_{\xi}=2$\ and $p=3$,
so then $n_{\xi}=10$ and $n_{\nu}=m_{\xi}+1=3$. We used a Smolyak sparse grid
with Gauss\textendash Legendre quadrature points for collocation.
For the solution of the Navier\textendash Stokes problem we used the hybrid strategy
with $20$ steps of Picard iteration followed by at most $20$ steps of Newton
iteration. We used mean viscosity $\nu_{1}=4.5455\times10^{-3}$, which
corresponds to Reynolds number $Re=220$, and the rightmost eigenvalue is
$5.7963\times10^{-4}$ (the second largest eigenvalue is $-8.2273\times10^{-2}%
$), see the right panel in Figure~\ref{fig:lambdaM}.
Table~\ref{tab:symstep} presents the results of validation and assessment of
the surrogates using the indicators from Section~\ref{sec:validation}. 
The trends are similar to those for 
the flow around an obstacle problem. For both $CoV$
$1\%$\ and $10\%$ the corresponding values of $\mu$ and $\sigma$ are in 
close agreement. The\ values of $\operatorname*{RMSE}$ are small for all
surrogates and again, they are smaller than the corresponding values of $\mu$
(and $\sigma$) by at least two orders of magnitude. Finally,
Figure~\ref{fig:symstep-pdf-Re} displays the probability density function
(pdf) estimates of the rightmost eigenvalue. We note that both pdf estimates
in this figure are ``narrower" comparing to the pdf estimates for flow around
an obstacle in Figure~\ref{fig:obstacle-pdf-Re}. Nevertheless there is an
excellent agreement of all estimates in both left and right panels.

\begin{table}[b]
\caption{Expansion flow around a symmetric step: validation of the surrogates
by Monte Carlo (MC) simulation using root mean square error (RMSE), their
assessment using estimates of the mean~$\mu$, standard deviation~$\sigma$, and
probability that the rightmost eigenvalue is nonnegative. The surrogates are
based on stochastic collocation (SC), Gaussian process regression (GP) and
neural network (NN), and the measures are defined in
Section~\ref{sec:validation}.}%
\label{tab:symstep}
\begin{center}
\begin{tabular}
[c]{|c|c|c|c|c|}\hline
& MC & SC & GP & NN\\\hline
\multicolumn{5}{|c|}{$CoV=1\%$}\\\hline
$\operatorname{RMSE}$ & - & $8.8129\times10^{-10}$ & $4.0545\times10^{-7}$ &
$1.5824\times10^{-8}$\\\hline
$\mu$ & $5.7982\times10^{-4}$ & $5.7982\times10^{-4}$ & $5.7987\times10^{-4}$
& $5.7982\times10^{-4}$\\
$\sigma$ & $2.9150\times10^{-4}$ & $2.9150\times10^{-4}$ & $2.9151\times
10^{-4}$ & $2.9149\times10^{-4}$\\\hline
$\operatorname{Pr}(\lambda\geq0)$ & \multicolumn{2}{|c|}{$98.4\%$} & $98.5\%$
& $98.4\%$\\\hline
\multicolumn{5}{|c|}{$CoV=10\%$}\\\hline
$\operatorname{RMSE}$ & - & $2.6183\times10^{-7}$ & $2.7106\times10^{-6}$ &
$4.2076\times10^{-7}$\\\hline
$\mu$ & $4.9677\times10^{-4}$ & $4.9676\times10^{-4}$ & $4.9711\times10^{-4}$
& $4.9685\times10^{-4}$\\
$\sigma$ & $2.9048\times10^{-3}$ & $2.9048\times10^{-3}$ & $2.9050\times
10^{-3}$ & $2.9048\times10^{-3}$\\\hline
$\operatorname{Pr}(\lambda\geq0)$ & \multicolumn{4}{|c|}{$57.5\%$}\\\hline
\end{tabular}
\end{center}
\end{table}

\paragraph{Computational time}

We briefly mention our experience with running the \textsc{MATLAB} functions
on a MacBook Pro laptop 
with a 3.5 GHz Intel Core i7 processor and 16 GB RAM.
The computation of the rightmost eigenvalue for one sample of $\xi$ using the simulator took at least 
$30%
\operatorname{s}%
$, depending on the value of$~\xi$ and settings of the inner solvers 
for the nonlinear iteration and call of the eigenvalue solver. 
On the other hand, a run of the emulators to evaluate the three surrogates took only 
between $0.02%
\operatorname{s}%
$ and $0.04%
\operatorname{s}%
$ for all $10^{3}$ sample points, which were used for validation and assessment. 
The learning part (construction of an emulator) took $0.18%
\operatorname{s}%
$ in case of~$\eta_{\operatorname*{GP}}$\ using the function \texttt{fitrgp}, 
and $1.25%
\operatorname{s}%
$ in case of~$\eta_{\operatorname*{NN}}$\ using the function
\texttt{train}. 
The construction of $\eta_{\operatorname*{SC}}$ was implemented as a part
of the simulator, however it can be seen, comparing~(\ref{eq:quadrature})
to~(\ref{eq:gPC-surrogate}), that if $n_{q}\approx n_{\xi}$ 
the construction of~$\eta_{\operatorname*{SC}}$ is inexpensive,
and in particular the timings of the construction of~$\eta_{\operatorname*{SC}}$ and its use are similar. 
Finally, we note that all three emulators were trained using only $n_q=29$ samples that require run of the simulator.
Therefore, since the overhead associated with the training and use of the emulators is very small, 
the computational savings provided by the emulators are dramatic.

\subsection{Effect of larger stochastic dimensions}

We also studied an effect of reducing the number of training (or design) points for the Gaussian process (GP) regression and neural network (NN) surrogates
using a problem with increasing stochastic dimension. 
We do not drop any quadrature (collocation) points from the stochastic collocation (SC) method since it would yield an incorrect quadrature rule. 
In particular, we considered the flow around an obstacle problem in a similar setup as in Section~\ref{sec:obstacle}
except with a channel of length~$12$ (instead of~$8$, cf. Figure~\ref{fig:mesh-obstacle}). 
There are then $12,640$ velocity and $1640$ pressure degrees of freedom, 
and the rightmost eigenvalue corresponding to the problem with the mean viscosity is a pair $0.0090\pm 2.2550i$. 
We considered a sequence of stochastic dimensions~$m_\xi=2,3,4,5$. 
Sizes of the gPC bases and numbers of the quadrature points are given in Table~\ref{tab:quadrature}. 
Other settings were the same as in Section~\ref{sec:obstacle}. 
We selected a fraction of the quadrature points to train the two surrogates for each of the stochastic dimensions 
in order to test the robustness in training of the Gaussian process regression and neural network surrogates.
For example, we selected every $10${\it th} quadrature point to be included in the training set, so that then the ratio $n_t/n_q =10\%$.  
Tables~\ref{tab:obstacle-2} and~\ref{tab:obstacle-5} summarize the results for $m_\xi=2$ and $m_\xi=5$, respectively. 
From Table~\ref{tab:obstacle-2} it can be seen that by using only~$6$ training points, i.e., reducing the ratio $n_t/n_q$ to approximately $20\%$,
the GP surrogate already provides relatively a quite accurate estimate as compared to the results of the Monte Carlo simulation, 
whereas the results of the NN surrogate are not satisfactory. By increasing the number of the training points to $8$ leads to a dramatic 
improvements of the NN surrogate. Nevertheless, by including all quadrature points into the training set, 
the approximation provided by the NN appears to be slightly more accurate then the one provided by the GP regression,
but overall the most accurate is the result provided by the stochastic collocation. 
The same trends can be observed also from Table~\ref{tab:obstacle-5} for the case with $m_\xi=5$,
except that in this case only approximately $5\%$ of the quadrature points are needed for the GP regression to provide a reasonable surrogate,
and approximately $10\%$ are needed for the NN. Therefore it appears that either of the GP or NN surrogates may provide an attractive alternative 
to the stochastic collocation for the high-dimensional problems.

\begin{table}[b]
\caption{Sizes of the gPC bases~$n_\xi$ and numbers of the quadrature points~$n_q$ for stochastic dimensions~$m_\xi$ and gPC degree~$p=3$.}%
\label{tab:quadrature}
\begin{center}
\begin{tabular}
[c]{|c|c|c|c|c|c|}\hline
$m_\xi$ & $1$ & $2$ & $3$ & $4$ & $5$ \\ \hline
$n_\xi$ & $4$ & $10$ & $20$ & $35$ & $56$ \\ \hline
$n_q$ & $4$ & $29$ & $69$ & $137$ & $241$ \\ \hline
\end{tabular}
\end{center}
\end{table}

\begin{table}[b]
\caption{Effect of reducing the number of training points on the GP and NN surrogates 
for the flow around an obstacle problem with the channel of length $12$ and with $m_\xi=2$.
The same quantities are used as in Table~\ref{tab:obstacle}, and they were defined in Section~\ref{sec:validation}.}%
\label{tab:obstacle-2}
\begin{center}
\begin{tabular}
[c]{|c|c|c|c|c|}\hline
\multirow{3}{*}{MC} & $\mu$  & \multicolumn{3}{|c|}{$8.8125 \times 10^{-3}$}  \\ \cline{2-5}
   & $\sigma$ & \multicolumn{3}{|c|}{$7.1136 \times 10^{-3}$} \\ \cline{2-5}
   & $\operatorname{Pr}(\lambda\geq0)$ & \multicolumn{3}{|c|}{$89.7\%$} \\ \hline
\multirow{5}{*}{GP} & $n_d \,\, (\approx n_d/n_q)$ & $6 \, (\approx 20\%)$ & $8 \, (\approx 30\%)$ & $29 \, (100\%)$ \\ \cline{2-5}
   & $\operatorname{RMSE}$ & $5.7386 \times 10^{-6}$ & $3.8676 \times 10^{-6}$ & $2.4706 \times 10^{-6}$ \\ \cline{2-5}
   & $\mu$ & $8.8134 \times 10^{-3}$ & $8.8103 \times 10^{-3}$ & $8.8117 \times 10^{-3}$ \\ \cline{2-5}
   & $\sigma$ & $7.1124 \times 10^{-3}$ & $7.1135  \times 10^{-3}$ & $7.1135 \times 10^{-3}$ \\ \cline{2-5}
   & $\operatorname{Pr}(\lambda\geq0)$  & \multicolumn{3}{|c|}{$89.7\%$}  \\ \hline
\multirow{5}{*}{NN} & $n_t \,\, (\approx n_t/n_q)$ & $6 \, (\approx 20\%)$ & $8 \, (\approx 30\%)$ & $29 \, (100\%)$ \\ \cline{2-5}
   & $\operatorname{RMSE}$ & $6.1469 \times 10^{-3}$ & $7.7102  \times 10^{-5}$ & $1.4824 \times 10^{-7}$ \\ \cline{2-5}
   & $\mu$ & $13.6704 \times 10^{-3}$ & $8.8029 \times 10^{-3}$ & $8.8125 \times 10^{-3}$ \\ \cline{2-5}
   & $\sigma$ & $4.8058 \times 10^{-3}$ & $7.1469 \times 10^{-3}$ & $7.1135 \times 10^{-3}$ \\ \cline{2-5}
   & $\operatorname{Pr}(\lambda\geq0)$  & $100\%$ & $89.4\%$ & $89.7\%$ \\ \hline
\multirow{5}{*}{SC} & $n_q$ & \multicolumn{3}{|c|}{$29$} \\ \cline{2-5}
   & $\operatorname{RMSE}$ & \multicolumn{3}{|c|}{$3.0072 \times 10^{-8}$} \\ \cline{2-5}
   & $\mu$ & \multicolumn{3}{|c|}{$8.8125 \times 10^{-3}$} \\ \cline{2-5}
   & $\sigma$ & \multicolumn{3}{|c|}{$7.1136 \times 10^{-3}$} \\ \cline{2-5}
   & $\operatorname{Pr}(\lambda\geq0)$ & \multicolumn{3}{|c|}{$89.7\%$} \\ \hline
\end{tabular}
\end{center}
\end{table}

\begin{table}[b]
\caption{Effect of reducing the number of training points on the GP and NN surrogates 
for the flow around an obstacle problem with the channel of length $12$ and with $m_\xi=5$.
The same quantities are used as in Table~\ref{tab:obstacle}, and they were defined in Section~\ref{sec:validation}.}%
\label{tab:obstacle-5}
\begin{center}
\begin{tabular}
[c]{|c|c|c|c|c|}\hline
\multirow{3}{*}{MC} & $\mu$  & \multicolumn{3}{|c|}{$8.7886 \times 10^{-3}$}   \\ \cline{2-5}
   & $\sigma$ & \multicolumn{3}{|c|}{$9.3753 \times 10^{-3}$} \\ \cline{2-5}
   & $\operatorname{Pr}(\lambda\geq0)$ & \multicolumn{3}{|c|}{$82.4\%$}  \\ \hline
\multirow{5}{*}{GP} & $n_d \,\, (\approx n_d/n_q)$ & $13 \, (\approx 5\%)$ & $25 \, (\approx 10\%)$ & $241 \, (100\%)$ \\ \cline{2-5}
   & $\operatorname{RMSE}$ & $1.7554 \times 10^{-3}$ & $ 1.2636 \times 10^{-5}$ & $1.2535 \times 10^{-5}$  \\ \cline{2-5}
   & $\mu$ & $8.8913 \times 10^{-3}$ & $8.7898 \times 10^{-3}$ & $8.7914 \times 10^{-3}$  \\ \cline{2-5}
   & $\sigma$ & $9.1899 \times 10^{-3}$ & $9.3740 \times 10^{-3}$ & $9.3742 \times 10^{-3}$  \\ \cline{2-5}
   & $\operatorname{Pr}(\lambda\geq0)$  & $83.0\%$ & \multicolumn{2}{|c|}{$82.4\%$}  \\ \hline
\multirow{5}{*}{NN} & $n_t \,\, (\approx n_t/n_q)$ & $13 \, (\approx 5\%)$ & $25 \, (\approx 10\%)$ & $241 \, (100\%)$ \\ \cline{2-5}
   & $\operatorname{RMSE}$ & $7.6909 \times 10^{-3}$ & $9.3273 \times 10^{-5}$ & $5.1869 \times 10^{-6}$ \\ \cline{2-5}
   & $\mu$ & $9.3639 \times 10^{-3}$ & $8.7848 \times 10^{-3}$ & $8.7886  \times 10^{-3}$ \\ \cline{2-5}
   & $\sigma$ & $1.7564 \times 10^{-3}$ & $9.3215  \times 10^{-3}$ & $9.3717 \times 10^{-3}$ \\ \cline{2-5}
   & $\operatorname{Pr}(\lambda\geq0)$ & $100\%$ & $82.3\%$ & $82.4\%$\\ \hline
\multirow{5}{*}{SC} & $n_q$ & \multicolumn{3}{|c|}{$241$} \\ \cline{2-5}
   & $\operatorname{RMSE}$ & \multicolumn{3}{|c|}{$1.7987 \times 10^{-7}$} \\ \cline{2-5}
   & $\mu$ & \multicolumn{3}{|c|}{$8.7886 \times 10^{-3}$} \\ \cline{2-5}
   & $\sigma$ & \multicolumn{3}{|c|}{$9.3754 \times 10^{-3}$} \\ \cline{2-5}
   & $\operatorname{Pr}(\lambda\geq0)$ & \multicolumn{3}{|c|}{$82.4\%$} \\ \hline
\end{tabular}
\end{center}
\end{table}

\section{Conclusion}

\label{sec:conclusion}We studied linear stability of Navier\textendash Stokes
equations with stochastic viscosity. This leads to a generalized eigenvalue
problem, and we are interested in characterization of the rightmost
eigenvalue. 
We designed three emulators for construction of the
rightmost eigenvalue surrogate. The first surrogate was based on generalized
polynomial chaos, and it was constructed using stochastic collocation, resp.
its pseudospectral variant (sometimes called nonintrusive stochastic Galerkin
method), which uses integration on Smolyak sparse grid and numerical
quadrature. For the second and third surrogates we used functions available in
\textsc{Matlab}. The second surrogate was based on Gaussian process
regression, and we used function \texttt{fitrgp}. The third surrogate was
based on shallow neural network, and we used function \texttt{fitnet} with
Bayesian Regularization backpropagation. We found that the set of quadrature
points used for the generalized polynomial chaos surrogate is also suitable
for training the other two emulators (based on Gaussian processes and neural
network), and we also found that certain scaling of the learning data points,
and subsequent descaling of the predictions, proposed by these emulators,
improves the quality of the surrogates. Finally, for the benchmark problems,
all three surrogates were in excellent agreement with the Monte Carlo simulation,
and we also found that the number of training points used for 
the Gaussian process regression and the neural network can be further reduced 
without compromising the quality of the surrogates.


{\small
{\em Authors' addresses}: \\
{\em Bed\v{r}ich Soused\'{i}k}, University of Maryland, Baltimore County, MD, USA \\
 e-mail: \texttt{sousedik@\allowbreak umbc.edu}. \\
{\em Howard C. Elman}, University of Maryland, College Park, MD, USA \\
 e-mail: \texttt{helman@\allowbreak umd.edu}. \\
{\em Kookjin Lee}, Arizona State University, Tempe, AZ, USA \\
 e-mail: \texttt{kookjin.lee@\allowbreak asu.edu}. \\
 {\em Randy Price}, George Mason University, Fairfax, VA, USA \\
 e-mail: \texttt{rprice25@\allowbreak gmu.edu}.
}

\end{document}